\magnification=1100 
\font\bigbf=cmbx10 scaled \magstep2
\font\medbf=cmbx10 scaled \magstep1 
\hfuzz=10 pt
\baselineskip=16pt
\input psfig.sty

\font\medbf=cmbx10 at 13pt

\def\i{\item}
\def\n{\noindent}

\def\L{{\bf L}}
\def\meas{{\rm meas}}

\def\forall{\hbox{~~for all}~~}
\def\sqr#1#2{\vbox{\hrule height .#2pt
\hbox{\vrule width .#2pt height #1pt \kern #1pt
\vrule width .#2pt}\hrule height .#2pt }}
\def\square{\sqr74}
\def\endproof{\hphantom{MM}\hfill\llap{$\square$}\goodbreak}
\def\ve{\varepsilon}
\def\n{\noindent}
\def\c{\centerline}

\def\F{{\cal F}}

\def\I{{\cal I}}
\def\O{{\cal O}}
\def\S{{\cal S}}
\def\wto{\rightharpoonup}
\def\Ups{\Upsilon}
\def\sgn{{\rm sgn}}
\def\C{{\cal C}}
\def\R{I\!\!R}
\def\implies{\Longrightarrow}

\def\meas{\hbox{meas}\,}

\def\tv{\hbox{Tot.Var.}}
\def\vs{\vskip 2em}
\def\vsk{\vskip 4em}
\def\v{\vskip 1em}
\null 
\c{\bigbf A Sharp Decay Estimate}
\v
\c{\bigbf
for Positive Nonlinear Waves}
\vs 
\c{\it Alberto Bressan$^{(*)}$ and
Tong Yang$^{(**)}$}
\v
\c{(*) ~S.I.S.S.A., Via Beirut 4, Trieste 34014, ITALY}
\v
\c{(**)~ Department of Mathematics, City University of Hong Kong}
\vsk
\n{\bf Abstract.}  
We consider a strictly hyperbolic, genuinely nonlinear system
of conservation laws in one space dimension. A sharp decay estimate is
proved for the positive waves in an entropy weak solution.
The result is stated in terms of a partial ordering among
positive measures, using
symmetric rearrangements and a comparison with a solution of
Burgers' equation with impulsive sources.
\vsk 
\n{\medbf 1 - Introduction}
\v
Consider a strictly hyperbolic system of $n$ conservation laws
$$u_t+f(u)_x=0\eqno(1.1)$$
and assume that all characteristic fields are genuinely nonlinear.
Call $\lambda_1(u)<\cdots<\lambda_n(u)$ the eigenvalues
of the Jacobian matrix $A(u)\doteq Df(u)$. 
We shall use bases of left and right eigenvectors $l_i(u)$, $r_i(u)$
normalized so that
$$\nabla\lambda_i(u)\, r_i(u)\equiv 1\,,\qquad\qquad 
l_i(u)\,r_j(u)=\cases{1\quad &if\quad $i=j$,\cr
0\quad &if\quad $i\not=j$.\cr}\eqno(1.2)$$
Given a function $u:\R\mapsto\R^n$ with small total variation, 
following [BC], [B] one can define 
the measures
$\mu^i$ of $i$-waves in $u$ as follows.
Since $u\in BV$, its distributional derivative $D_xu$ is a Radon measure.
We define $\mu^i$ as the measure such that
$$\mu^i\doteq l_i(u)\cdot D_xu\eqno(1.3)$$
restricted to the set where $u$ is continuous, while, at each point $x$ 
where $u$ has a jump, we define
$$\mu^i\big(\{x\}\big)\doteq\sigma_i\,,\eqno(1.4)$$
where $\sigma_i$ is the strength of the $i$-wave in the solution of
the Riemann problem with data 
$u^-=u(x-)$, $u^+=u(x+)$. 
In accordance with (1.2), if the solution of the Riemann problem contains 
the intermediate states $u^-=\omega_0, \omega_1,\ldots,\omega_n=u^+$, 
the strength of the $i$-wave is defined as 
$$\sigma_i\doteq \lambda_i(\omega_i)-\lambda_i(\omega_{i-1}).\eqno(1.5)$$ 
Observing that 
$$\sigma_i= l_i(u^+)\cdot (u^+ - u^- )+ O(1)\cdot 
|u^+ - u^- |^2, 
$$ 
we can find a vector $l_i(x)$ 
such that 
$$\left|l_i(x)-l_i\big( u(x+)\big)\right| 
=\O(1)\cdot \big| u(x+)-u(x-)\big|, 
\eqno(1.6)$$ 
$$\sigma_i= 
l_i(x)\cdot \big( u(x+)-u(x-)\big). 
\eqno(1.7)$$ 
We can thus define the measure $\mu^i$ equivalently as  
$$\mu^i\doteq l_i\cdot D_x u\,,\eqno(1.8)$$ 
where $l_i(x)=l_i\big( u(x)\big)$ at points where $u$ is 
continuous, while $l_i(x)$ is some 
vector which satisfies (1.6)-(1.7) at points of jump. 
For all $x\in\R$ 
there holds 
$$\Big|l_i(x)- l_i\big( u(x)\big)\Big|=\O(1)\cdot \big| 
u(x+)-u(x-)\big|\,.\eqno(1.9)$$

We call $\mu^{i+}$, $\mu^{i-}$ respectively  the positive  
and negative parts of $\mu^i$, so that 
$$\mu^i=\mu^{i+}-\mu^{i-},\qquad\qquad |\mu^i|=\mu^{i+}+\mu^{i-}.\eqno(1.10)$$ 
It is our purpose to prove a sharp estimate on the decay of the density  
of the measures $\mu^{i+}$. 
This will be achieved by introducing a partial ordering within the  
family of positive Radon measures. In the following, $meas(A)$ denotes the 
Lebesgue measure of a set $A$. 
\v 
\n{\bf Definition 1.} {\it  
Let $\mu,\mu'$ be two positive Radon measures. 
We say that 
$\mu\preceq \mu'$ if and only if  
$$\sup_{meas(A)\leq s}\mu(A)~\leq~\sup_{meas(B)\leq s}\mu'(B) 
\qquad\qquad\hbox{for every}~s>0\,. 
\eqno(1.11)$$ 
} 
 
In some  sense, the 
above relation means that $\mu'$ is more singular than $\mu$.  
Namely, it has a greater total mass, concentrated   
on regions with higher density. 
Notice that the usual order relation 
$$\mu\leq\mu'\qquad\hbox{if and only if}\qquad \mu(A)\leq \mu'(A) \quad 
\hbox{for every}\quad A\subset\R\eqno(1.12)$$ 
is much stronger. Of course 
$\mu\leq\mu'$ implies $\mu\preceq\mu'$, 
but the converse does not hold. 
 
Following [BC], [B], together with the measures $\mu^i$ we define the 
Glimm functionals 
$$V(u)\doteq \sum_i |\mu^i|(\R)\,,\eqno(1.13)$$ 
$$Q(u)\doteq \sum_{i<j} \big(|\mu^j|\otimes|\mu^i|\big) \big\{ 
(x,y)\,;~x<y\big\} 
+\sum_i\big(\mu^{i-}\otimes|\mu^i|\big) \big\{ 
(x,y)\,;~x\not=y\big\}\,. 
\eqno(1.14)$$ 
Let now $u=u(t,x)$ be an entropy weak solution 
of (1.1).  If the total variation of $u$ is small and the constant 
$C_0$ is large enough, 
it is well known that the quantities 
$$Q(t)\doteq Q\big(u(t)\big)\,,\qquad\qquad 
\Ups(t)\doteq V\big(u(t)\big)+C_0\,Q\big(u(t)\big)\eqno(1.15)$$ 
are non-increasing in time.    
The decrease in $Q$ controls the amount of interaction, while 
the decrease in $\Ups$ controls both the 
interaction and the cancellation in the solution. 
 
An accurate estimate on the measure $\mu^{i+}_t$ of positive $i$-waves in  
$u(t,\cdot)$ will be obtained by a comparison with a solution 
of Burgers' equation with source terms. 
\v 
\n{\bf Theorem 1.} {\it  For some constant 
$\kappa$ and for every small BV solution $u=u(t,x)$ 
of the system (1.1) the following holds. 
Let $w=w(t,x)$ be the solution of the scalar Cauchy problem 
with impulsive source term 
$$w_t+(w^2/2)_x=-\kappa\,\sgn(x)\cdot {d\over dt}Q\big(u(t)\big)\,,
\eqno(1.16)$$ 
$$w(0,x)=\sgn(x)\cdot \sup_{meas(A)<2|x|} 
{\mu^{i+}_0(A)\over 2}\,.\eqno(1.17)$$ 
Then, for every $t\geq 0$,  
$$\mu^{i+}_t\preceq D_x w(t)\,.\eqno(1.18)$$ 
} 
\v 
As shown in the 
next section, the initial data in (1.17) represents the  
{\it odd rearrangement} of the 
function $v_i(x)\doteq \mu^{i+}_0\big(]-\infty,x]\big)$. 
The above theorem improves the earlier estimate derived in [BC].  
For a scalar conservation law with strictly convex flux,  
a classical decay estimate was proved by Oleinik [O]. 
In the case of genuinely nonlinear systems, results related to the decay 
of nonlinear waves were also obtained in [GL], [L1], [L2], [BG]. 
An application of the present analysis will appear in [BY], 
where Theorem 1 is used to estimate the rate of convergence 
of vanishing viscosity approximations. 
\vsk 
\n{\medbf 2 - Lower semicontinuity} 
\v 
Let $\mu$ be a positive Radon measure on $\R$, so that 
$\mu \doteq D_x v$ is the distributional derivative of some 
bounded, non-decreasing function 
$v:\R\mapsto\R$. 
We can decompose 
$$\mu=\mu^{\rm sing}+\mu^{ac}$$ 
as the sum of a singular and an absolutely continuous part, w.r.t.~Lebesgue 
measure.  The absolutely continuous part corresponds to 
the usual derivative $z\doteq v_x$, which is a non-negative $\L^1$ function 
defined at a.e.~point. We shall denote by $\hat z$ the {\it 
symmetric rearrangement} of $z$, i.e.~the unique even function 
such that 
$$\hat z(x)=\hat z(-x)\,,\qquad\qquad \hat z(x)\geq \hat z(x')~~~ 
\hbox{if}~~~0<x<x'\,, 
\eqno(2.1)$$ 
$$ 
\meas\Big(\big\{x\,;~\hat z(x)> c\big\}\Big)= 
\meas\Big(\big\{x\,;~z(x)> c\big\}\Big)\qquad\qquad \hbox{for every}~ c>0 
\,. 
\eqno(2.2)$$ 
Moreover, we define the {\it odd rearrangement} of $v$ as the  
unique function 
$\hat v$ such that (fig.~1) 
$$\hat v(-x)=-\hat v(x)\,,\qquad\qquad 
\hat v(0+) = {1\over 2}\mu^{\rm sing}(\R)\,,\eqno(2.3)$$ 
$$\hat v(x)=\hat v(0+)+\int_0^x z(y)\,dy \qquad 
\hbox{for}~x>0\,.\eqno(2.4)$$ 
By construction, the function $\hat v$ is convex for $x<0$ and  
concave for $x>0$.   
 
\midinsert 
\vskip 10pt 
\centerline{\hbox{\psfig{figure=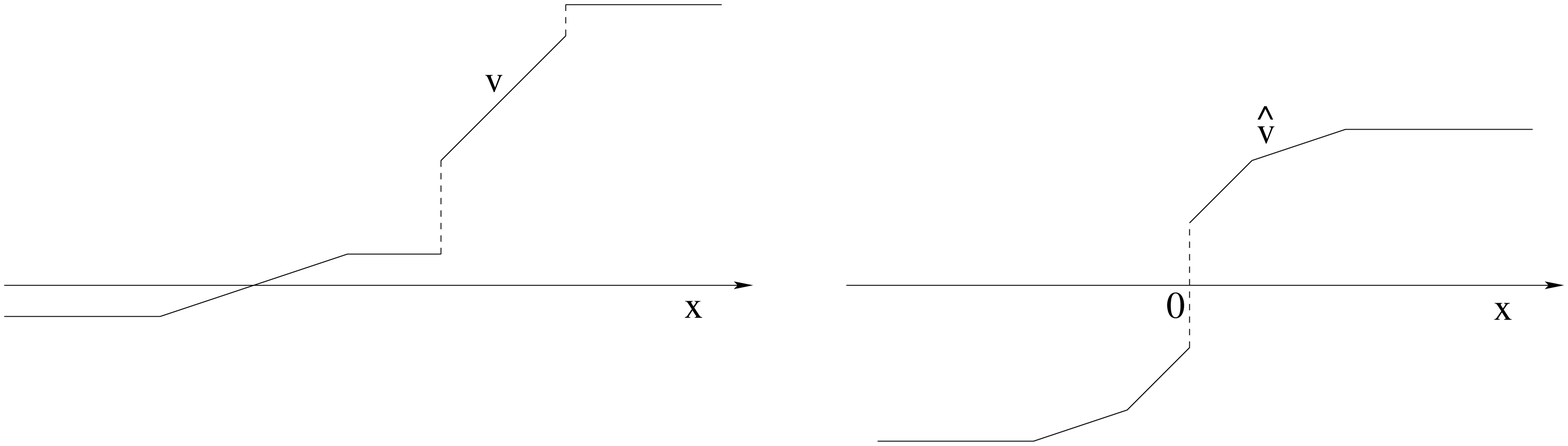,width=14cm}}} 
\centerline{\hbox{figure 1}} 
\vskip 10pt 
\endinsert 
 
The relation between the odd rearrangement $\hat v$ 
and the partial ordering (1.10) is clarified by the following 
result, which is an easy consequence of the definitions. 
\v 
\n{\bf Proposition 1.} {\it 
Let $\mu=D_x v$ and 
$\mu'=D_x v'$ be  positive Radon measures. 
Call $\hat v,\hat v'$ the odd rearrangements
of $v,v'$, respectively.   Then 
$\mu\preceq D_x \hat v\preceq \mu$ and moreover
$$\hat v(x)=\sgn(x)\cdot\sup_{meas (A)\leq 2|x|} 
{\mu(A)\over 2}\,,
\eqno(2.5)$$ 
$$\mu\preceq\mu'\qquad\hbox{if and only if}\qquad
\hat v(x)\leq \hat v'(x)\quad\hbox{for all}~x>0\,.\eqno(2.6)$$
} 
\v 
Two more results will be used in the sequel. By the restriction 
of a measure $\mu$ to a set $J$, we mean the measure 
$$(\mu\lfloor J) (A)\doteq \mu(A\cap J)\,.$$ 
\v 
\n{\bf Proposition 2.} {\it Let $\mu,\mu'$ be  
positive measures. Consider any finite partition 
$\R=J_1\cup\cdots\cup J_N$.  If the restrictions of $\mu,\mu'$ 
to each set $J_\ell$ satisfy 
$\mu\lfloor{J_\ell}~\preceq~ \mu'\lfloor{J_\ell} 
\,$, 
then $\mu\preceq\mu'$.} 
\v 
\n{\bf Proposition 3.} {\it Assume that $\mu\preceq D_s w$ for some 
nondecreasing odd function $w$.  If 
$|\mu^\sharp-\mu|(\R)\leq \ve$, then 
$$\mu^\sharp\preceq D_s\left[ w+\sgn(s)\cdot {\ve\over 2}\right]\,.$$ 
} 
\v 
The next result is concerned with  
the lower semicontinuity of the partial ordering 
$~\preceq~$ w.r.t.~weak convergence of measures. 
\v 
\n{\bf Proposition 4.}  
{\it Consider a sequence of 
measures $\mu_\nu$ converging weakly to a measure $\mu$. 
Assume that the positive parts satisfy  
$\mu_\nu^+\preceq D w_\nu$ for some odd, nondecreasing 
functions $s\mapsto w_\nu(s)$, concave for $s>0$.  Let $w$ be the odd funtion 
such that 
$$w(s)\doteq\liminf_{\nu\to\infty} w_\nu(s)\qquad\qquad\hbox{for}~~ s>0\,.$$ 
Then the positive part of $\mu$ satisfies 
$$\mu^+\preceq D_s w\,.\eqno(2.7)$$ 
}\v 
\n{\bf Proof.} By possibly taking a subsequence, we can assume that 
$w_\nu(s)\to w(s)$ for all $s\not= 0$.  Moreover, we can assume the weak 
convergence 
$$\mu_\nu^+\wto \tilde \mu^+\,,\qquad\qquad 
\mu_\nu^-\wto \tilde \mu^-\,,$$ 
for some positive measures $\tilde\mu^+$,  $\tilde\mu^-$. 
We thus have 
$$\mu=\tilde\mu^+-\tilde\mu^-\,,\qquad\qquad 
\mu^+\leq\tilde \mu^+\,,\qquad  \mu^-\leq\tilde\mu^-\,.\eqno(2.8)$$ 
By (2.8) it suffices to prove that $\tilde\mu^+\preceq D_s w$, i.e. 
$$\hbox{meas}\,(A)\leq 2s\qquad\implies\qquad 
\tilde \mu^+(A)\leq 2w(s)\,,\eqno(2.9)$$ 
for every $s>0$ and every Borel measurable set $A\subset\R$. 
If (2.9) fails, there exists $s>0$ and a set $A$ 
such that 
$$\meas(A)=2s\,,\qquad\tilde\mu^+(A)>2w(s)=2\lim_{\nu\to\infty}w_\nu(s). 
$$ 
Since $w$ is continuous for $s>0$, 
we can choose an open set $A'\supseteq A$ such that, setting 
$s'\doteq\meas(A')/2$, one has 
$2w(s')<\tilde \mu^+(A)$.   
By the weak convergence $\mu^+_\nu\rightharpoonup \tilde\mu^+$ one obtains 
$$\tilde \mu^+(A')\leq\liminf_{\nu\to\infty} \mu_\nu^+(A')\leq 
2w(s') <\tilde \mu^+(A)\,,$$ 
reaching a contradiction. Hence (2.9) must hold. 
\endproof 
\v 
Toward the proof of Theorem 1 we shall need 
a lower semicontinuity 
property for wave measures, similar to what proved in [BaB].  
In the following, $C_0$ is the same constant as in (1.15).  
\v 
\n{\bf Lemma 1.}  {\it  
Consider a sequence of functions  $u_\nu$ 
with uniformly small total variation and call 
$\mu^{i+}_\nu$ the corresponding measures of positive 
$i$-waves.  Let $s\mapsto w_\nu(s)$, 
$\nu\geq 1$, be a sequence of odd, nondecreasing functions, 
concave for $s>0$, such that  
$$\mu^{i+}_\nu\preceq D_s \Big[ w_\nu +C_0\, \sgn(s)  
\big(Q_0-Q(u_\nu)\big)\Big]\,\eqno(2.10)$$ 
for some $Q_0$. 
Assume that $u_\nu\to u$ and $w_\nu\to w$ in $\L^1_{\rm loc}$. 
Then the measure of positive $i$-waves in $u$ satisfies 
$$\mu^{i+}\preceq D_s \Big[ w +C_0\, \sgn(s)  
\big(Q_0-Q(u)\big)\Big]\,.\eqno(2.11)$$ 
} 
\v 
\n{\bf Proof.}   
The main steps follow the proof of Theorem 10.1 in [B]. 
\v 
\n{\bf 1.}  
By possibly taking a 
subsequence we can assume that $u_\nu(x)\to u(x)$ for  
every $x$  
and that the measures of total variation 
converge weakly, say 
$$|\mu_\nu|\doteq 
\big|D_x u_\nu\big|\wto\mu^\sharp\eqno(2.12)$$  
for some positive Radon measure $\mu^\sharp$.  In this 
case one has $\mu^\sharp\geq |\mu|$, in the sense of (1.12). 
\v 
\n{\bf 2.} 
Let any 
$\ve>0$ be given.  
Since the total mass of $\mu^\sharp$ is finite, one can select 
finitely many points $y_1,\dots,y_N$ such that 
$$\mu^\sharp\big(\{x\}\big)<\ve,\qquad\quad \forall~ x\notin\{ 
y_1,\ldots, y_N\}.\eqno(2.13)$$ 
We now choose disjoint open intervals $I_k\doteq 
\,]y_k-\rho,~y_k+\rho[\,$ 
such that 
$$\mu^\sharp\big(I_k\setminus\{y_k\}\big) <{\ve\over 
N}\qquad\qquad k=1,\ldots,N.\eqno(2.14)$$ 
Moreover, we choose $R>0$ such that 
$$\bigcup_{k=1}^N I_k\subset [-R,R], 
\qquad\quad\mu^\sharp\big(]-\infty,\,-R]\cup [R,\,\infty[\big)  
< \ve. \eqno(2.15) 
$$ 
Because of (2.13), we can now 
choose points 
$p_0<-R<p_1<\cdots<R<p_r$ which 
are continuity points for $u$ and for every $u_\nu$,  
such that  
$$\mu^\sharp\big(\{p_h\}\big)=0,\quad\qquad 
u_\nu(p_h)\to u(p_h)\qquad \forall h=0,\ldots,r, 
\eqno(2.16)$$ 
and such that either 
$$p_h-p_{h-1}<{\ve\over N}\,,\qquad\quad 
p_{h-1}<y_k<p_h,\qquad\quad [p_{h-1},p_h]\subset I_k,\eqno(2.17)$$ 
for some $k\in\{1,\ldots,N\}$, or else 
$$|\mu|\big([p_{h-1},p_h]\big)\leq 
\mu^\sharp\big([p_{h-1},p_h]\big)<\ve.\eqno(2.18)$$ 
Call $J_h\doteq [p_{h-1},p_h]$. 
If (2.18) holds, by weak convergence 
for some $\nu_0$ sufficiently large one has 
$$|\mu_\nu| (J_h) < \ve \qquad\qquad \forall\nu\geq\nu_0.\eqno(2.19)$$ 
On the other hand, if (2.17) holds, from (2.14) it follows 
$$|\mu| \big(J_h\setminus \{y_k\}\big)\leq 
\mu^\sharp\big(J_h\setminus\{y_k\}\big) < {\ve\over N}\,.\eqno(2.20)$$ 
\v 
In the remainder of the proof, the main strategy is as follows. 
 
\i{$\bullet$} On the intervals $J_{h(k)}$ containing a point $y_k$  
of large oscillation, we first replace each $u_\nu$ by a piecewise constant 
function $\bar u_\nu$ having a single jump at $y_k$. 
The relations between the corresponding measures 
$\mu^i_\nu$ and $\bar \mu^i_\nu$ are given by Lemma 10.2 in [B]. 
Then we take the limit as $\nu\to\infty$. 
 
\i{$\bullet$} On the remaining intervals $J_h$ with small oscillation, 
we replace the left eigenvectors $l_i(u_\nu)$ by a constant vector 
$\l_i(u_h^*)$. Then we use Proposition 4 to estimate the limit as  
$\nu\to\infty$. 
\v 
\n{\bf 3.}  
We first take care of the intervals $J_h$ containing a point $y_k$ of 
large oscillation, so that (2.17) holds.   
For each $k=1,\ldots,N$, let $h=h(k)\in\{1,\ldots,r\}$ 
be the index such that 
$y_k\in J_h\doteq [p_{h-1},~p_h]$. 
For every $\nu\geq 1$ consider the function 
$$\bar u_\nu(x)\doteq\cases{u_\nu(x) &if~~$x\notin\cup_k J_{h(k)}$, 
\cr 
u_\nu(p_{h(k )-1}) &if~~$x\in \,]p_{h(k)-1},~y_k[\,$, \cr 
u_\nu(p_h) &if~~$x\in\,[y_k,~p_{h(k)}]$. \cr}$$ 
Observe that all functions $u, \bar u_\nu$ are continuous at every point 
$p_0,\ldots,p_r$ and have jumps at $y_1,\ldots, y_N$.   
Call $\bar\mu_\nu^{i}$, $i=1,\ldots,n$, the 
corresponding measures, defined as in (1.8) with 
$u$ replaced by $\bar u_\nu$.  Clearly $\bar\mu_\nu^{i}=\mu_\nu^i$ 
outside the intervals $J_{h(k)}$ of large oscillation. 
By Lemma 10.2 at p.203 in [B], there holds 
$$Q(\bar u_\nu)\leq Q(u_\nu),\qquad\qquad 
V(\bar u_\nu)+C_0\, Q(\bar u_\nu)\leq 
V(u_\nu)+C_0 \cdot Q(u_\nu),$$ 
$$\bar \mu^{i+}_\nu(\R)- \mu^{i+}_\nu(\R)~\leq~ 
C_0\, \big[ Q(u_\nu)-Q(\bar u_\nu)\big].$$ 
As a consequence, from (2.10) we deduce 
$$\bar \mu^{i+}_\nu\preceq D_s \Big[ T^\ve w_\nu +C_0\, \sgn(s)  
\big(Q_0-Q(\bar u_\nu)\big)\Big]\,,\eqno(2.21)$$ 
where  
$$T^\ve w(s)\doteq\cases{w(s+\ve/2)\qquad &if\quad $s>0$,\cr 
w(s-\ve/2)\qquad &if\quad $s<0$.\cr}$$ 
Indeed, all the mass which in $\mu^{i+}_\nu$ lies 
on the set 
$$\Omega\doteq \bigcup_{k=1}^N J_{h(k)}\,,\qquad\qquad 
J_h\doteq [p_{h-1},~p_h]$$ 
is replaced in $\bar\mu^{i+}_\nu$ by 
point masses at $y_1,\ldots,y_N$.  
We obtain (2.21) by observing that, by (2.17), 
$meas(\Omega)<\ve$. Moreover, the increase in the total mass 
is $\leq C_0\big[Q(u_\nu)-Q(\bar u_\nu)\big]$. 
 
Since $u_\nu(p_h)\to u(p_h)$ for every $h$, there holds 
$$\eqalign{\Big| \mu^i\big(\{ y_k\}\big)-\bar\mu_\nu^i\big(\{y_k\} 
\big)\Big| 
&=\O(1)\cdot\Big\{ \big| u(y_k-)-u(p_{h(k)-1})\big| 
+\big| u(y_k+)-u(p_{h(k)})\big|\cr 
& \qquad\quad +\big| u(p_{h(k)-1})-u_\nu(p_{h(k)-1})\big| 
+\big| u(p_{h(k)})-u_\nu(p_{h(k)})\big|\Big\}\cr 
&=\O(1)\cdot{\ve\over N}\cr}\eqno(2.22)$$ 
for each $k=1,\ldots,N$ and all $\nu$ sufficiently large. 
By construction we also have 
$$|\bar \mu_\nu^i|\big( J_{h(k)}\setminus\{ y_k\}\big)=0, 
\qquad\quad |\mu^i|\big( J_{h(k)}\setminus\{ y_k\}\big)=\O(1)\cdot {\ve 
\over N}\,.\eqno(2.23)$$ 
\v 
\n{\bf 4.}  
Next, call $\S\doteq\big\{h\,;~\mu^\sharp(J_h)<\ve\big\}$  
the family of intervals where the oscillation of every $u_\nu$ is small, 
so that (2.18) holds. 
If $h\in\S$, for every $x,y\in J_h$ and $\nu$ sufficiently large, 
one has 
$$\big|u_\nu(x) - u_\nu(y)\big|\leq |\mu_\nu | (J_h) 
<\ve, 
$$ 
$$\big|u(x) - u(y)\big|\leq |\mu | (J_h)\leq \mu^\sharp(J_h) 
<\ve.  
$$ 
Set $u^*_h\doteq 
u(p_h)$. By the pointwise convergence $u_\nu(p_h)\to u(p_h)$ 
and the two above estimates it follows 
$$\big|u_\nu(x) - u_h^* \big| <\ve,\qquad 
\big| u(x)-u^*_h\big| <\ve,\qquad\quad \forall x\in J_h\,. 
\eqno(2.24) 
$$ 
\v 
\n{\bf 5.} 
We now introduce the measures $\hat \mu^i_\nu$ 
such that 
$$\hat\mu^i_\nu\doteq l_i(u_h^*)\cdot D_x u_\nu$$ 
restricted to each interval $J_h$, $h\in\S$ where the oscillation  
is small, while 
$$\hat\mu^i_\nu=\bar\mu^i_\nu$$ 
on each interval $J_h=J_{h(k)}$ where the oscillation is large. 
Observe that the restriction of $\hat \mu^i_\nu$ to $J_{h(k)}$ 
consists of a single mass at the point $y_k$.  
Namely,  
$\hat \mu^i_\nu\big(\{y_k\}\big)$ is precisely the size of the 
$i$-th wave in the solution of the Riemann problem 
with data $u^-=u_\nu(p_{h(k)-1})$,~$u^+=u_\nu(p_{h(k)})$. 
 
We define $\hat w_\nu$ as the non-decreasing odd function such that 
$$\hat w_\nu(s)~\doteq\sup_{meas(A)\leq 2s}~ {\hat \mu^{i+}_\nu(A)\over 2} 
\,,\qquad\qquad s>0.\eqno(2.25)$$ 
By possibly taking a further subsequence we can assume 
the convergence 
$$Q(\bar u_\nu)\to \overline Q\,,\qquad\qquad 
\hat \mu^i_\nu\wto \hat \mu^i\,,\qquad\qquad \hat w_\nu(s)\to \hat w(s)\,. 
$$ 
Using (2.16), we can apply Proposition 4 on each interval $J_h$ 
and obtain 
$$\hat\mu^{i+}\preceq D_s\hat w\,.\eqno(2.26)$$ 
\v 
\n{\bf 6.} 
Observe that, by (2.24) and (2.19), 
$$|\hat \mu^i_\nu-\mu^i_\nu|(J_h)=\O(1)\cdot \ve\,\mu^\sharp(J_h)\qquad\qquad 
h\in\S\,,\eqno(2.27)$$ 
{}From (2.21) and the definition of $\hat w_\nu$ at (2.25) 
it thus follows 
$$\hat w_\nu(s)\leq T^\ve w_\nu(s)+C_0\big[Q_0-Q(\bar u_\nu)\big]+ 
\O(1)\cdot\ve 
\qquad\qquad s>0\,.\eqno(2.28)$$ 
Letting $\nu\to\infty$ we obtain 
$$\hat w(s)\leq T^\ve w(s)+C_0[Q_0-\overline Q]+ 
\O(1)\cdot\ve\qquad\qquad s>0\,, 
\eqno(2.29)$$ 
$$\overline Q=\lim_{\nu\to\infty} Q(\bar u_\nu)\geq 
\lim_{\nu\to\infty} Q(u_\nu)-\O(1)\cdot\ve \geq 
Q(u)-\O(1)\cdot\ve\,,\eqno(2.30)$$ 
because of the lower semicontinuity of the functional $u\mapsto Q(u)$. 
{}From (2.26), (2.29) and (2.30) we deduce 
$$\hat\mu^{i+}\preceq D_s\Big[ T^\ve w+\sgn (s)\,\big( C_0[Q_0-Q(u)] 
+\O(1)\cdot\ve\big)\Big]\,.$$ 
By (2.22)--(2.24), our construction of the measure $\hat \mu^i$  
achieves the property 
$$\big|\mu^{i+}-\hat\mu^{i+}\big|(\R)=\O(1)\cdot\ve\,.$$ 
Hence, by Proposition 3, 
$$\mu^{i+}\preceq D_s\Big[ T^\ve w+\sgn (s)\,\big( C_0[Q_0-Q(u)] 
+\O(1)\cdot\ve\big)\Big]\,.$$ 
Since $\ve>0$ was arbitrary, this proves (2.11). 
\endproof 
\vsk 
\n{\medbf 3 - A decay estimate} 
\v 
The second basic ingredient in the proof  
is the following lemma, which 
refines the estimate in [BC]. 
\v 
\n{\bf Lemma 2.} {\it For some constant $\kappa>0$ the following holds.  
Let $u=u(t,x)$ be any entropy weak solution 
of (1.1), with initial data $u(0,x)=\bar u(x)$ having small total variation. 
Then the measure $\mu^{i+}_t$ 
of positive $i$-waves in $u(t,\cdot)$ can be estimated as follows. 
 
Let $w:[0,\tau[\,\times\R\mapsto\R$ be the solution of 
Burgers' equation  
$$w_t+(w^2/2)_x=0\eqno(3.1)$$ 
with initial data 
$$w(0,x)=\sgn (x)\cdot 
\sup_{meas(A)\leq 2|x|}\,{\mu^{i+}_0(A)\over 2}\,.\eqno(3.2)$$ 
Set 
$$w(\tau,x)=w(\tau-\,,~x)+ 
\kappa\,\sgn(x)\cdot \big[ Q(\bar u)-Q\big(u(\tau)\big)\Big]\,.\eqno(3.3)$$ 
Then  
$$\mu^{i+}_\tau\preceq D_x w(\tau)\,.\eqno(3.4)$$ 
} 
\v 
\n{\bf Proof.} 
The main steps follow the proof of Theorem 10.3 in [B]. 
We first prove the estimate (3.3) under the additional hypothesis: 
\v 
\item{(H)} There exist points  
$y_1<\cdots <y_m$ such that the initial data 
$\bar u$ is smooth outside such points, constant for  
$x<y_1$ and $x>y_m$, and 
the derivative component $l_i(u)\,u_x$ is constant on each 
interval $\,]y_{\ell},\,y_{\ell+1}[\,$.   
Moreover, the Glimm functional $t\mapsto Q\big(u(t)\big)$ 
is continuous at $t=\tau$. 
\v 
\n{\bf 1.}  
The solution $u=u(t,x)$ can be obtained as limit of  
front tracking approximations. In particular, 
we can consider a particular converging  
sequence $(u_\nu)_{\nu\geq 1}$ of $\ve_\nu$-approximate solutions 
with the following additional properties: 
 
\i{(i)} Each $i$-rarefaction front $x_\alpha$ travels with the  
characteristic speed of the state on the right: 
$$\dot x_\alpha=\lambda_i\big( u(x_\alpha+)\big).$$ 
 
\i{(ii)} Each $i$-shock front $x_\alpha$ travels with a speed strictly 
contained between the right and the left characteristic speeds: 
$$\lambda_i\big( u(x_\alpha+)\big)<\dot x_\alpha< 
\lambda_i\big( u(x_\alpha-)\big).\eqno(3.5)$$ 
 
\i{(iii)} As $\nu\to\infty$, the interaction potentials satisfy 
$$Q\big(u_\nu(0,\cdot)\big) 
\to Q(\bar u).\eqno(3.6)$$ 
\v 
\n{\bf 2.} Let $u_\nu$ be an approximate 
solution constructed by the front tracking algorithm. 
By a {\it (generalized) $i$-characteristic}  
we mean an absolutely continuous curve 
$x=x(t)$ such that 
$$\dot x (t)\in \big[ \lambda_i(u_\nu(t,x-)),~ 
\lambda_i(u_\nu(t,x+))\big]$$ 
for a.e.~$t$.  If $u_\nu$ satisfies the above properties (i)-(ii), then  
the $i$-characteristics are precisely the 
polygonal lines $x:[0,\tau]\mapsto \R$ for which the following holds. 
For a suitable partition $0=t_0<t_1<\cdots <t_m=\tau$, on each subinterval 
$[t_{j-1},\,t_j]$ either $\dot x(t)=\lambda_i\big(u_\nu(t,x)\big)$, or 
else $x$ coincides with a wave-front of the $i$-th family. 
For a given terminal point $\bar x$ we shall consider the 
{\it minimal backward $i$-characteristic} through $\bar x$, defined as 
$$y(t)=\min\big\{ x(t)\,; ~~x~\hbox{ is an $i$-characteristic}, 
~~x(\tau)=\bar x\big\}.$$ 
Observe that $y(\cdot)$ is itself an $i$-characteristic. 
By (3.5), it cannot coincide with an $i$-shock front of $u$  
on any nontrivial time interval. 
 
In connection with the exact solution $u$, we define an $i$-characteristic 
as a curve  
$$t\mapsto x(t)=\lim_{\nu\to\infty}x_\nu(t)$$  
which is the limit of 
$i$-characteristics in a sequence of front tracking solutions 
$u_\nu\to u$. 
\v 
\n{\bf 3.}  
Let $\ve>0$ be given.  If the assumption (H) holds, 
the measure $\mu^{i+}_\tau$ of $i$-waves 
in $u(\tau)$ is supported on a bounded interval and 
is absolutely continuous w.r.t.~Lebesgue 
measure.   We can thus find a piecewise constant function 
$\psi^\tau$ with jumps at points $x_1(\tau)<\bar x_2(\tau)<\ldots<\bar  
x_N(\tau)$ such that 
$$\int \left| {d\mu^{i+}_\tau\over dx}-\psi^\tau\right|\,dx <\ve\,,\qquad 
\qquad \int_{x_j(\tau)} 
^{x_{j+1}(\tau)} 
\left( {d\mu^{i+}_\tau\over dx}-\psi^\tau\right)\,dx=0\qquad j=1,\ldots,N-1\,. 
\eqno(3.7)$$

\midinsert 
\vskip 10pt 
\centerline{\hbox{\psfig{figure=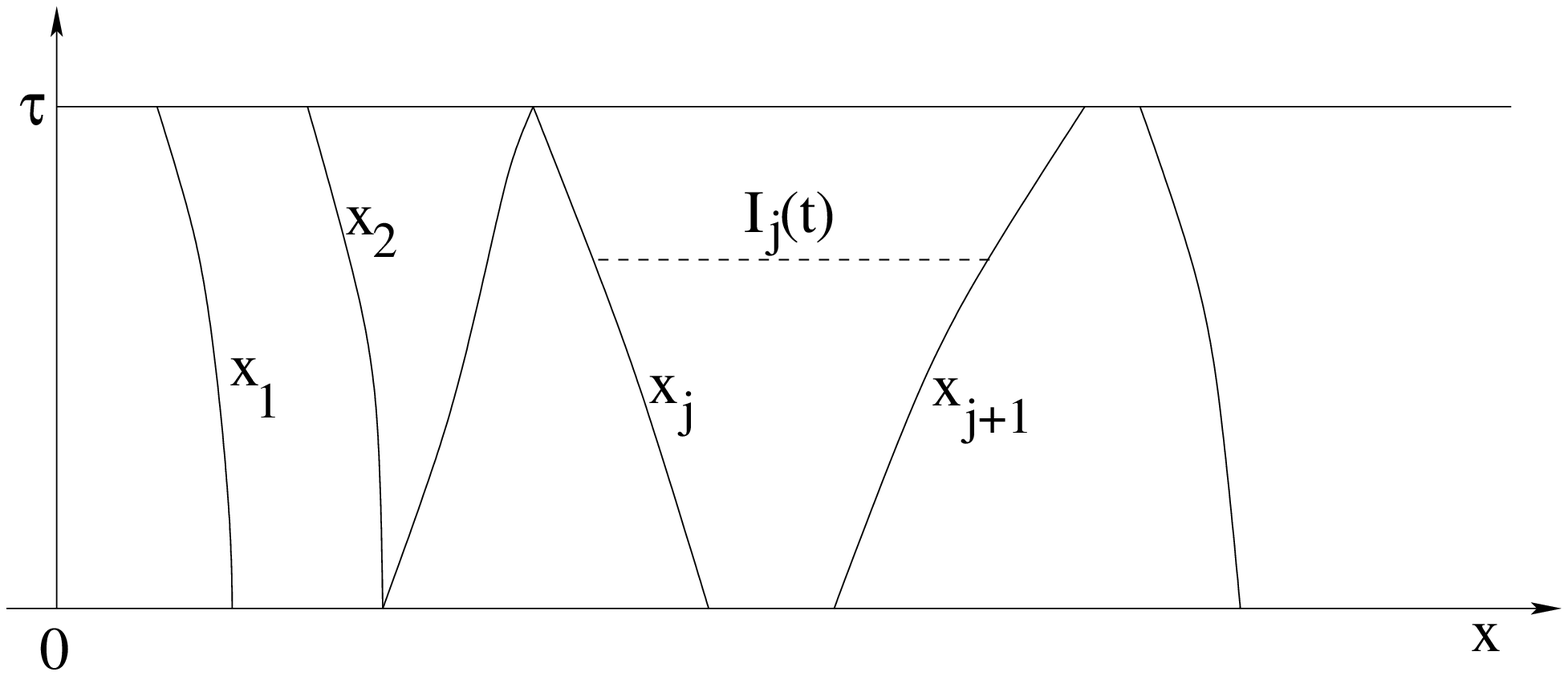,width=10cm}}} 
\centerline{\hbox{figure 2}} 
\vskip 10pt 
\endinsert 
\v

To prove the lemma in this special case, relying on 
Proposition 2, it thus suffices to 
find $i$-characteristics 
$t\mapsto x_j(t)$ such that the following holds (fig.~2) 
\v 
\i{(i)} 
For each $j=1,\ldots,N$, the function 
$\psi^\tau$ is constant on the interval  
$\,\big]x_j(\tau),\,x_{j+1}(\tau)\big[\,$ and (3.7) holds.  
Moreover, either 
$x_j(0)=x_{j+1}(0)$, or else the derivative component 
$\psi^0\doteq l^i(u) u_x(0,\cdot)$ is constant 
on the interval $\,\big]x_j(0),\,x_{j+1}(0)\big[\,$. 
\v 
\i{(ii)}  
An estimate corresponding to (3.3)-(3.4) holds restricted 
to each subinterval $\big[ x_j(\tau),~x_{j+1}(\tau)\big[$. 
\v 
We need to explain in more detail  
this last statement. 
Define 
$$I_j(t)\doteq \big[ x_j(t),\,x_{j+1}(t)\big[\,,\qquad\qquad 
\Delta_j\doteq\big\{ (t,x)\,;~~t\in [0,\tau]\,, 
~~x\in I_j(t)\big\}\,.$$ 
For each $j$, we denote by $\Gamma_j$ the total amount 
of wave interaction within the domain $\Delta_j$. This is defined 
as in [B], first for a sequence of front tracking approximations 
$u_\nu$, 
then taking a limit as $\nu\to\infty$. 
Furthermore, we define the constant values 
$$\psi_j^\tau\doteq \psi^\tau(x)\qquad\qquad x\in I_j(\tau)\,,$$ 
$$\psi_j^0\doteq \psi^0(x)\qquad\qquad x\in I_j(0)\,,$$ 
Call $$\sigma_j^0\doteq \lim_{t\to 0+}\mu^{i+}\big(I_j(t)\big)$$ 
the initial amount of positive $i$-waves 
inside the interval $I_j$. 
 
For each interval $I_j$, we consider on one hand  
the function $w^\tau_j$  
corresponding to (3.2)-(3.3), namely 
$$w^\tau_j(s)\doteq \min\left\{ \sigma_j^0\,,~ 
{s\over \tau+(\psi_j^0)^{-1}}\right\} 
+\kappa\Gamma_j\cdot \sgn(s)\,.$$ 
Here $(\psi_j^0)^{-1}\doteq 0$ in the case where 
$x_j(0)=x_{j+1}(0)$.  This may happen when the initial data has a jump 
at $x_j(0)$, and the corresponding measure $\mu^{i+}$ has  
a Dirac mass (with infinite density) at that point. 
 
On the other hand, we look at the nondecreasing, odd function 
$\eta_j$ such that 
$$\eta_j(s)\doteq \min\Big\{ \psi^\tau_j\,s,~~\psi_j^\tau\big[x_{j+1}(\tau)- 
x_j(\tau)\big]\Big\}\qquad\qquad s>0\,.$$ 
Our basic goal is to prove that (fig.~3) 
$$\eta_j(s)\leq w^\tau_j(s)\qquad\qquad\hbox{for all}~s>0\,.\eqno(3.8)$$ 
Indeed, by (3.7), for $s>0$ one has  
$$ 
\sup_{meas (A)\leq 2s} {\mu^{i+}_\tau\big(A\cap I_j(\tau)\big)\over 2} 
\leq \eta_j(s)+\ve_j$$ 
with $$\sum_j \ve_j<\ve\,.$$ 
Proving (3.8) for each $j$ will thus imply 
$$ 
\mu^{i+}_\tau \preceq w(\tau,x)=w(\tau-\,,~x)+ 
\kappa\,\sgn(x)\cdot \big[ Q(\bar u)-Q\big(u(\tau)\big)+\O(1)\cdot\ve\Big] 
\,.$$ 
Since $\ve>0$ was arbitrary, this establishes the lemma under the 
additional assumptions (H). 
 
\midinsert 
\vskip 10pt 
\centerline{\hbox{\psfig{figure=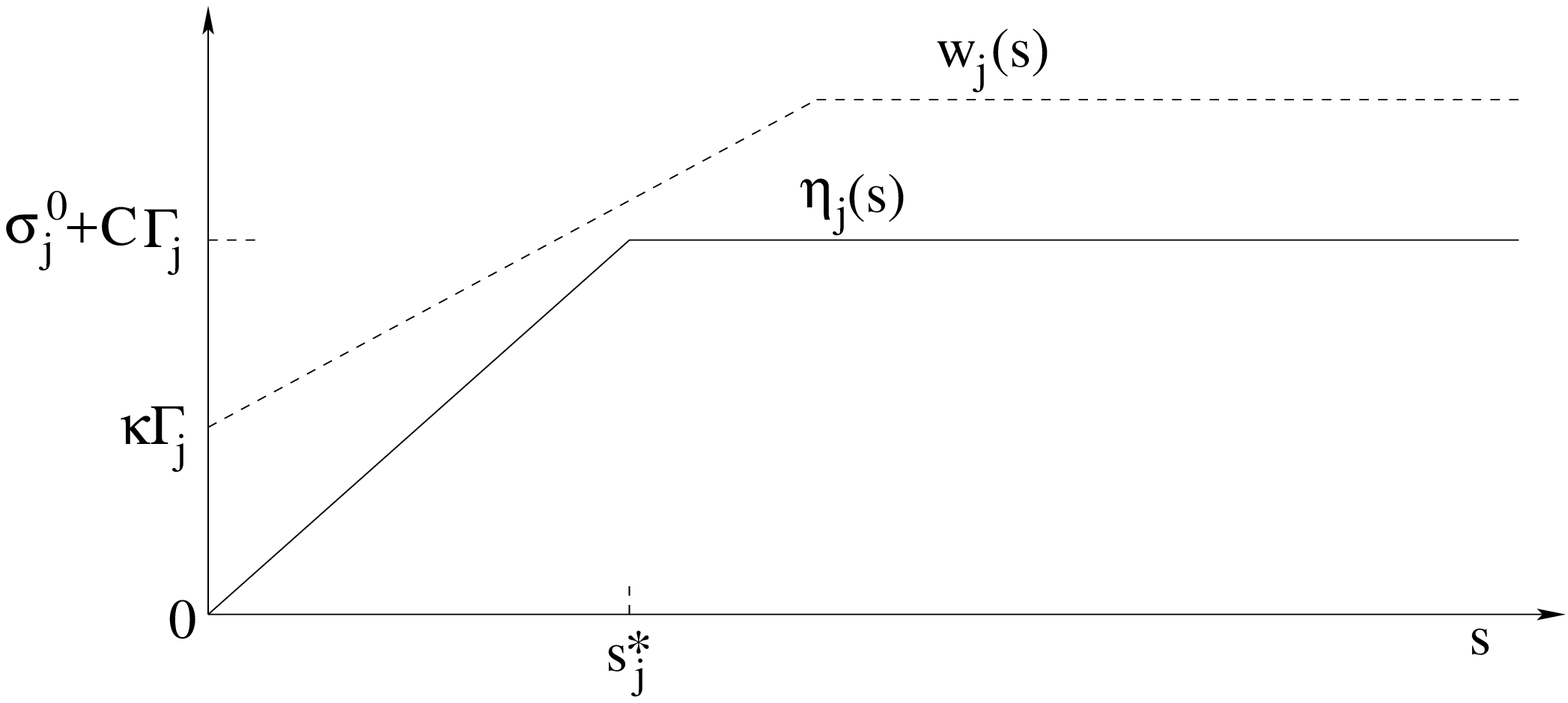,width=12cm}}} 
\centerline{\hbox{figure 3}} 
\vskip 10pt 
\endinsert 
 
\v 
\n{\bf 4.} We now work toward a proof of (3.8), in three cases. 
\v 
\n Case 1:~ 
$\sigma_j^0=0$. 
\v 
\n Case 2:~ $x_j(0)=x_{j+1}(0)$ and $\sigma_j^0>0$. 
\v 
\n Case 3:~ $x_j(0)<x_{j+1}(0)$ and $\sigma_j^0 
=\big(x_{j+1}(0)-x_j(0)\big)\,\psi_j^0>0$. 
\v 
In Case 1 the proof is easy. Indeed, the total amount 
of positive $i$-waves in $I_j(\tau)$ is here bounded by a constant times 
the total amount of interaction taking place inside the domain $\Delta_j$, 
i.e. 
$$\mu^{i+}_\tau\big(I_j(\tau)\big)\leq C_0\cdot \Gamma_j$$ 
for some constant $C_0$. On the other hand 
$$w^\tau_j(s)=\kappa\Gamma_j\cdot\sgn(s)\,.$$ 
Choosing $\kappa>C_0$ we achieve (3.8). 
\v 
\n{\bf 5.} 
Since Case 2 can be obtained from Case 3 in the limit  
as $x_{j+1}-x_j\to 0$, 
we shall only give a proof for Case 3.

We can again distinguish two cases.  If the amount 
of interaction  $\Gamma_j$ 
is large compared with the initial amount of $i$-waves, say  
$$\Gamma_j\geq {1\over 6C_0} \sigma_j^0\,,$$ 
then the bound (3.8) is readily achieved choosing $\kappa>8C_0$. 
Indeed, for $s>0$ we have  
$$\eta_j(s)\leq {1\over 2}\,\mu^{i+}_\tau\big(I_j(\tau)\big) 
\leq C_0\Gamma_j+\sigma_j^0\leq 7C_0\,\Gamma_j\,.$$ 
 
The more difficult case to analyse is when $\Gamma_j$ is 
small, say 
$$\Gamma_j<\sigma_j^0/6C_0\,.\eqno(3.9)$$ 
Looking at figure 3, it clearly suffices to prove (3.8) for the single value 
$$s=s_j^*\doteq {x_{j+1}(\tau)-x_j(\tau)\over 2}\,.$$ 
Equivalently, calling 
$$z_j(t)\doteq x_{j+1}(t)-x_j(t)$$ 
the length of the interval $I_j(t)$ and 
$$\sigma_j^\tau\doteq \mu^{i+}_\tau\big(I_j(\tau)\big)=z_j(\tau)\,\psi^\tau_j 
$$ 
the total amount of positive $i$-waves inside $I_j(\tau)$, 
we need to show that 
$$\sigma_j^\tau\leq 
2\kappa\Gamma_j+ \min\left\{ \sigma_j^0\,,~~{2s_j^*\over 
\tau +(\psi_j^0)^{-1}}\right\}\,.\eqno(3.10)$$ 
 
By the approximate conservation of $i$-waves over the region $\Delta_j$, we 
can write 
$$\sigma_j^\tau\leq \sigma_j^0+C_0\Gamma_j\,.\eqno(3.11)$$ 
Using (3.11) in (3.10), our task is reduced to showing that 
$$\sigma_j^\tau\leq 
2\kappa\Gamma_j+ {2s_j^*\over 
\tau +(\psi_j^0)^{-1}}\eqno(3.12)$$ 
for a suitably large constant $\kappa$. 
Because of (3.11), it suffices to show that 
$$\eqalign{z_j(\tau)&\geq (\sigma_j^0-C'\Gamma_j)\big(\tau+(\psi_j^0)^{-1}\big) 
\cr 
&=\big[ z_j(0)+\tau\sigma_j^0\big] -C'\big(\tau+(\psi_j^0)^{-1}\big)\, 
\Gamma_j 
\cr}\eqno(3.13)$$ 
for a suitable constant $C'$. 
 
\midinsert 
\vskip 10pt 
\centerline{\hbox{\psfig{figure=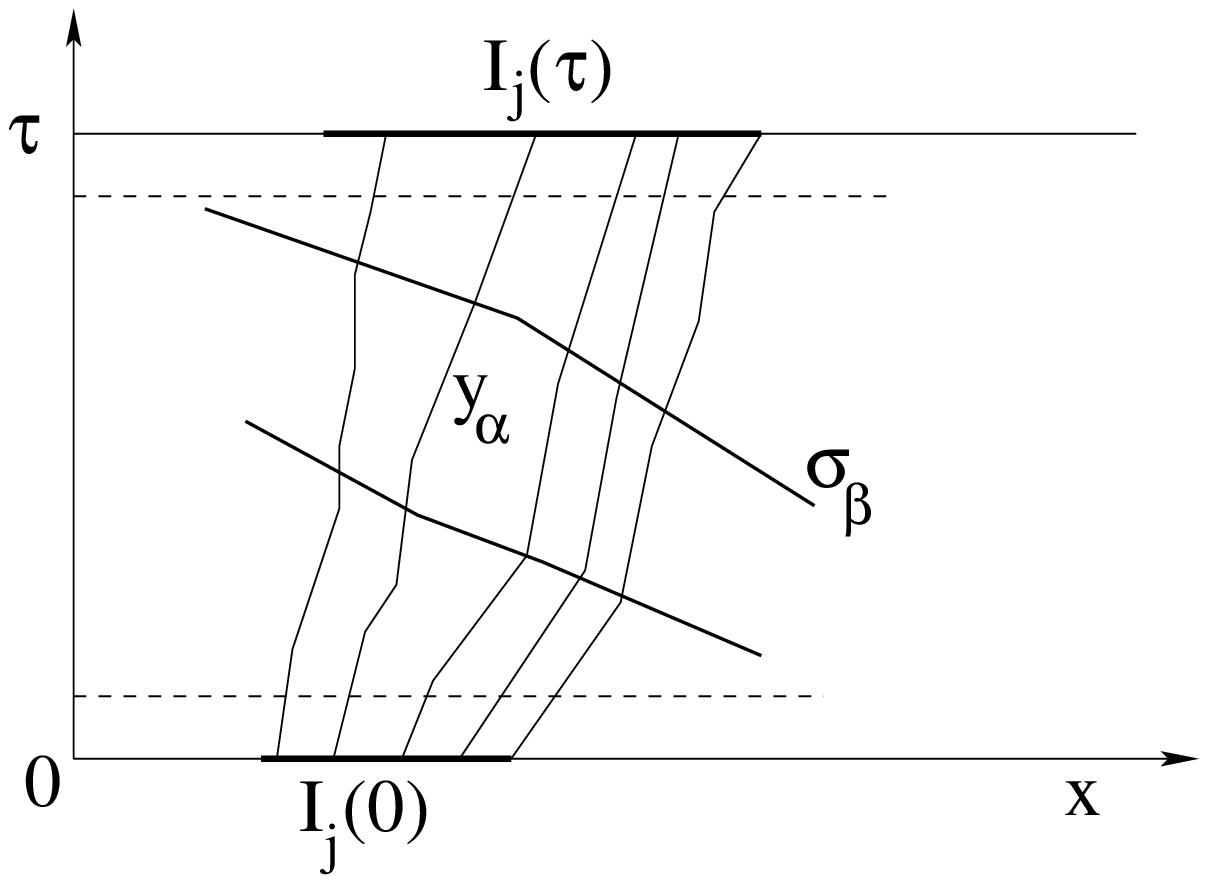,width=10cm}}} 
\centerline{\hbox{figure 4}} 
\vskip 10pt 
\endinsert 
 
\v 
\n{\bf 6.} We now prove (3.13). 
Notice that, by genuine nonlinearity and the normalization (1.2), 
if no other waves were present in the region $\Delta_j$ 
we would have $\Gamma_j=0$ and 
$${d\over dt}z_j(t)\equiv \sigma_j^0\,.$$ 
In this case, the equality would hold in (3.13). 
 
To handle the general case, 
we represent the solution $u$ as a limit of front tracking 
approximations $u_\nu$, where for each $\nu\geq 1$ the function 
$u_\nu(0,\cdot)$ contains exactly $\nu$ rarefaction fronts 
equally spaced along the interval $I_j(0)$. 
Each of these fronts has initial strength 
$\sigma_\alpha(0)=\sigma^0_j/\nu$. 
For $\alpha=1,\ldots,\nu,$ let $y_\alpha(t)\in I_j(t)$ be  
the location of one of these fronts at time $t\in [0,\tau]$, 
and let $\sigma_\alpha(t)>0$ be its strength. 
Moreover, call 
$$J_\alpha(t)\doteq \big[ y_\alpha(t)\,,~y_{\alpha+1}(t)\big[\,, 
\qquad\qquad \Delta_\alpha\doteq \big\{ (t,x)\,;~~t\in [0,\tau]\,,~~ 
x\in J_\alpha(t)\big\}\,,$$ 
and let $\Gamma_\alpha$ be the total amount of interaction 
in $u_\nu$ taking place inside the domain $\Delta_\alpha$. 
  
We define a subset of indices $\I\subseteq 
\{1,\ldots,\nu\}$ 
by setting  $\alpha\in\I$ if 
$$5C_0\Gamma_\alpha>\sigma_\alpha(0)=\sigma_j^0/\nu\,.\eqno(3.14)$$ 
Observe that, if $\alpha\notin\I$, then 
$$\left|{\sigma_\alpha(t)\over  
\sigma_\alpha(0)}-1\right| <{1\over 2}\qquad 
\hbox{for all}~t\in [0,\tau]\,.$$   
In particular, if $\alpha,\alpha+1\notin\I$, then the interval $J_\alpha(t)$ 
is well defined for all $t\in [0,\tau]$. 
Its length  
$$z_\alpha(t)\doteq y_{\alpha+1}(t)-y_\alpha(t)$$ 
satisfies the differential inequality 
$${d\over dt}z_\alpha(t)\geq 
W_\alpha(t)-C_1\cdot \sum_{\beta\in\C_\alpha(t)}|\sigma_\beta|\eqno(3.15)$$ 
for some constant $C_1$. Here 
$$\eqalign{W_\alpha(t)&\doteq \big[\hbox{amount  
of $i$-waves inside the interval}~J_\alpha(t) 
\big]\cr 
&\geq \sigma_\alpha(0)-C_0\Gamma_\alpha\,,\cr}\eqno(3.16)$$ 
while $\C_\alpha(t)$ refers to the set of all  
wave fronts of different families which are crossing the interval $J_\alpha$ 
at time $t$. 
Calling $W_\alpha'$ the total amount of waves of families $\not= i$ 
which lie inside $J_\alpha(0)$, 
we can now write 
$$\int_0^\tau \left(\sum_{\beta\in\C_\alpha(t)}|\sigma_\beta|\right)\,dt 
\leq \left(\max_{t\in [0,\tau]}z_\alpha(t)\right)\cdot {2\nu\over\sigma_j^0} 
\cdot\Gamma_\alpha +\O(1)\cdot\tau \Gamma_\alpha+\O(1)\cdot 
\left({z_j(0)+1\over\nu} 
\right)\,W_\alpha'\,.\eqno(3.17)$$ 
Indeed, by strict hyperbolicity, every front $\sigma_\beta$  
of a different family 
can spend at most a time $\O(1)\cdot z_\alpha$ inside $J_\alpha$. 
Either it is located inside $J_\alpha$ already at time $t=0$, 
or else, when it enters, it crosses $y_\alpha$ or $y_{\alpha+1}$. 
In this case, since $\alpha,\,\alpha+1\notin\I$, by (3.14)  
it will produce an interaction of magnitude  
$|\sigma_\beta\,\sigma_\alpha|\geq |\sigma_\beta\cdot \sigma_j^0|/2\nu$. 
The second term on the right hand side of (3.17) takes care of the 
new wave fronts which are generated through interactions 
inside $J_\alpha$. 
The last term takes into account wave front of different 
families that initially lie already inside $J_\alpha$ at time $t=0$. 
Integrating (3.15) over the time interval $[0,\tau]$  
and using (3.16)-(3.17) one obtains 
$$z_\alpha(\tau)~\geq~ z_\alpha(0)+\tau{\sigma_j^0\over\nu}- 
\O(1)\cdot \tau\Gamma_\alpha-\O(1)\cdot 
\left(\max_{t\in [0,\tau]}z_\alpha(t)\right)\cdot {2\nu\over\sigma_j^0} 
\cdot\Gamma_\alpha -\O(1)\cdot 
\left({z_j(0)+1\over\nu} 
\right)\,W_\alpha'\,.\eqno(3.18)$$ 
\v 
\n{\bf 7.} 
To proceed in our analysis, we now show that 
$$\max_{t\in [0,\tau]}z_\alpha(t)\leq 2\,z_\alpha(\tau)\,.\eqno(3.19)$$ 
Indeed, let $\tau'\in [0,\tau]$ be the time where the maximum is 
attained.  
If our claim (3.19) does not hold, there would exists a first time 
$\tau''\in [\tau',\tau]$ such that 
$z_\alpha(\tau'')=z_\alpha(\tau')/2$. 
>From (3.15) and the assumption $W_\alpha(t)\geq 0$  
it follows 
$$\int_{\tau'}^{\tau''} C_1\sum_{\beta\in\C_\alpha(t)}|\sigma_\beta|\,dt 
\geq {z_\alpha(\tau')\over 2}\,.\eqno(3.20)$$ 
Using the smallness of the total variation, a contradiction is now 
obtained as follows. 
Call 
$$\Phi(t)\doteq C_0Q(t)+\sum_{k_\beta\not= i}  
\phi_{k_\beta}\big(t,x_\beta(t)\big)\, 
|\sigma_\beta|\,,$$ 
where the sum ranges over all fronts of  
strength $\sigma_\beta$ located at $x_\beta$, of a family $k_\beta\not= i$. 
The weight functions $\phi_j$ are defined as 
$$\phi_j(t,x)\doteq \left\{\eqalign{0\qquad &\qquad 
\hbox{if}\quad x>y_{\alpha+1}(t)\,, 
\cr 
{y_{\alpha+1}(t)-x\over y_{\alpha+1}(t)-y_\alpha(t)} 
\quad &\qquad \hbox{if}\quad x\in \big[ 
y_\alpha(t),~y_{\alpha+1}(t)\big]\,, 
\cr 
1 \qquad&\qquad\hbox{if}\quad x<y_\alpha(t)\,, 
\cr}\right.$$ 
in the case $j>i$, while 
$$\phi_j(t,x)\doteq \left\{\eqalign{1\qquad &\qquad 
\hbox{if}\quad x>y_{\alpha+1}(t)\,, 
\cr 
{x-y_\alpha(t)\over y_{\alpha+1}(t)-y_\alpha(t)} 
\quad &\qquad\hbox{if}\quad x\in \big[ 
y_\alpha(t),~y_{\alpha+1}(t)\big]\,, 
\cr 
0\qquad &\qquad \hbox{if}\quad x<y_\alpha(t)\,, 
\cr}\right.$$ 
in the case $j<i$. 
Because of the term $C_0Q(t)$, the functional 
$\Phi$ is non-increasing at times of interactions. 
Moreover, outside interaction times  
a computation entirely similar to the one at p.213 of [B] 
now yields 
$$-{d\over dt} \Phi(t)\geq \sum_{\beta\in\C_\alpha(t)}|\sigma_\beta|\cdot 
{c_0\over z(t)}\,,\eqno(3.21)$$ 
for some small constant $c_0>0$ related to the gap between different 
characteristic speeds. 
From 
(3.20) and (3.21) respectively we now deduce 
$$\int_{\tau'}^{\tau''}  
\sum_{\beta\in\C_\alpha(t)}|\sigma_\beta|\,dt\geq  
{z_\alpha(\tau')\over 2C_1}\,,$$ 
$$\int_{\tau'}^{\tau''}  
\sum_{\beta\in\C_\alpha(t)}|\sigma_\beta|\,dt\leq 
\int_{\tau'}^{\tau''} \left|{d\Phi(t)\over dt}\right|\cdot {z_\alpha(\tau') 
\over c_0}\,dt\leq{ \Phi(\tau')\over c_0}\, z_\alpha(\tau')\,.$$ 
Since $\Phi(t)=\O(1)\cdot\tv\big\{u(t)\big\}$, by the smallness of the  
total variation we can assume 
$\Phi(\tau')<2C_1/c_0$.  In this case, the two above inequalities 
yield a contradiction. 
\v 
\n{\bf 8.} Using (3.19), from (3.18) we obtain 
$$\eqalign{ z_j(\tau)&=\sum_{1\leq\alpha\leq\nu} z_\alpha(\tau)\geq  
\sum_{\alpha\notin\I} z_\alpha(\tau)\cr 
&\geq  
\sum_{\alpha\notin\I}\left\{ {z_\alpha(0)+\tau\sigma^0_j/\nu 
\over 1+C_2 (\nu/\sigma_j^0)\Gamma_\alpha} -\O(1)\cdot \tau\Gamma_j 
-\O(1)\cdot \left({z_j(0)+1\over\nu}\right)W_\alpha'\right\}\cr 
&\geq  
\sum_{\alpha\notin\I}\left(z_\alpha(0)+\tau{\sigma^0_j\over\nu}\right) 
\left( 1-C_2{\nu\over\sigma_j^0}\Gamma_\alpha\right)- 
\O(1)\cdot \tau\Gamma_j-\O(1)\cdot {z_j(0)+1\over\nu} 
\cr 
&\geq  
\sum_{\alpha\notin\I}\left(z_\alpha(0)+\tau{\sigma^0_j\over\nu}\right) 
-C_2{z_j(0)\over \sigma_j^0}\Gamma_j 
-\O(1)\cdot \tau\Gamma_j-\O(1)\cdot {z_j(0)+1\over\nu} 
\,.\cr}\eqno(3.22)$$ 
By (3.14) the cardinality of the set $\I$ satisfies 
$$\#\I\cdot{\sigma_j^0\over 5C_0\nu}\leq\sum_{\alpha\in\I}\Gamma_\alpha\leq 
\Gamma_j\,,$$ 
hence 
$${\#\I\over\nu}\leq {5C_0\over \sigma_j^0}\Gamma_j\,.$$  
In turn, this implies 
$$ 
\sum_{\alpha\notin\I}\left(z_\alpha(0)+\tau{\sigma^0_j\over\nu}\right) 
\geq \big(z_j(0)+\tau\sigma_j^0\big)\left(1-{\#\I\over\nu} 
\right)\geq \big(z_j(0)+\tau\sigma_j^0\big)-5C_0\Gamma_j{z_j(0)\over\sigma_j^0} 
\Gamma_j-5C_0\tau\Gamma_j\,.\eqno(3.23)$$ 
Using (3.23) in (3.22), observing that 
$${z_j(0)\over \sigma_j^0}={x_{j+1}(0)-x_j(0)\over\sigma_j^0}=(\psi_j^0)^{-1}. 
$$ 
and letting $\nu\to\infty$ we conclude 
$$z_j(\tau)\geq \big(z_j(0)+\tau\sigma_j^0\big) 
-\O(1) \cdot (\psi_j^0)^{-1}\Gamma_j-\O(1)\cdot \tau\Gamma_j\,.$$ 
This establishes (3.13),for a suitable constant $C'$. 
\v 
\n{\bf 9.} In the general case,  
without the assumptions (H), the lemma is proved by an approximation 
argument.   We construct a convergent  
sequence of initial data $\bar u_\nu\to \bar u$ which satisfy 
(H) and such that 
$$\bar u_\nu\to \bar u\,,\qquad Q(\bar u_\nu)\to Q(\bar u)\,,\qquad 
\big|\mu^{i+}_{\nu,0}-\mu^{i+}_0\big|\to 0\,.$$ 
Calling $w_\nu$ the solution of  
(3.1) with initial data 
$$w_\nu(0,x)=\sgn (x)\cdot 
\sup_{meas(A)\leq 2|x|}\,{\mu^{i+}_{\nu,0}(A)\over 2}\,,$$ 
by the previous analysis we have 
$$\mu_{\nu,\tau_\nu}^{i+}\preceq D_x\Big[ 
w_\nu(\tau_\nu -)+\sgn(x)\cdot \big[ Q(\bar u_\nu)-Q(u_\nu(\tau_\nu))\big] 
\Big]\,.$$ 
Observe that $w_\nu(\tau_\nu-)\to w(\tau-)$ in $\L^1_{\rm loc}$.  
Choosing $\kappa\geq C_0$, by the lower semicontinuity 
result stated in Lemma 1 we now conclude 
$$\mu^{i+}_\tau\preceq 
D_x\Big[ 
w(\tau -)+\kappa\,\sgn(x)\cdot \big[ Q(\bar u)-Q(u(\tau))\big] 
\Big]\,.$$ 
\endproof 
\vsk 
 
\n{\medbf 4 - Proof of the main theorem} 
\v 
Using the previous lemmas, we now give a proof of Theorem 1. 
For a given interval $[0,\tau]$, 
the solution of the impulsive Cauchy problem (1.17)-(1.18)  
can be obtaines as follows. 
Consider a partition 
$0=t_0<t_1<\cdots<t_N=\tau$. 
Construct an approximate solution by 
requiring that $w(0,x)=\hat v_i(x)$, 
$$w_t+(w^2/2)_x=0\eqno(4.1)$$ 
on each subinterval $[t_{k-1}, t_k[\,$, while 
$$  
w(t_k,x)=w(t_k-\,,~x)+ 
\kappa\,\sgn(x)\cdot \big[ Q(t_{k-1})-Q(t_k)\big]\,.\eqno(4.2)$$ 
We then consider a sequence of  
partitions $0=t_0^\nu<t_1^\nu<\cdots<t_{N_\nu}^\nu=\tau$, and 
the corresponding 
solutions $w_\nu$. If the mesh of the partitions approaches zero, i.e. 
$$\lim_{\nu\to\infty}\sup_k |t_k^\nu-t^\nu_{k-1}| =0\,,$$ 
then the approximate solutions $w_\nu$ converge to a unique limit, 
which yields the solution of (1.17)-(1.18).  
 
Call $\F$  the set of nondecreasing odd functions, concave for  
$x>0$. This set is  
positively invariant for 
the flow of Burgers' equation (4.1).  
Moreover, this flow is order preserving.  Namely, 
if $w,w'\in\F$ are solutions of (4.1) with initial data such that 
$w(0,x)\leq w'(0,x)$ for all $x>0$, then  
also 
$$w(t,x)\leq w'(t,x)\qquad\qquad \hbox{for all}~~t,x>0\,.$$ 
Equivalently, 
$$D_x w(0)\preceq D_x w'(0)\qquad\implies D_x w(t)\preceq D_x w'(t)$$ 
for every $t>0$. 
For each fixed $\nu$, we can  
apply Lemma 2 on each subinterval $[t^\nu_{k-1}, t^\nu_k]$ 
and obtain 
$$\mu^{i+}_{t_k^\nu}\preceq D_x w_\nu(t_k^\nu)\qquad\implies\qquad 
\mu^{i+}_{t_{k+1}^\nu}\preceq D_x w_\nu(t_{k+1}^\nu)\,.$$ 
By induction on $k$, this yields 
$$\mu^{i+}_\tau\preceq D_x w_\nu(\tau)\,,\eqno(4.3)$$ 
where $w_\nu$ is the approximate solution constructed according to 
(4.1)-(4.2).   
Letting $\nu\to\infty$ and using Lemma 1, we achieve a proof of Theorem 1. 
\endproof

\vsk

{\bf Acknowledgments:}
The first author was supported by the Italian M.I.U.R.,
within the research project \# 2002017219  "Equazioni iperboliche e
paraboliche non lineari".
The research of the second author was supported by
CityU Direct Allocation Grant \# 7100198.        
\vsk

\c{\medbf References} 
\v
\i{[BaB]} P.~Baiti and A.~Bressan,
Lower semicontinuity of weighted path length in BV,
in ``Geometrical Optics and Related Topics",
F.~Colombini and N.~Lerner Eds., Birkh\"auser (1997), 
31-58.
\v
\i{[B]} A.~Bressan,
{\it Hyperbolic Systems of Conservation Laws. The One Dimensional
Cauchy Problem}, Oxford University Press, 2000.
\v
\i{[BC]} A.~Bressan and R.~M.~Colombo, 
Decay of positive waves in
nonlinear systems of conservation laws, {\it Ann. Scuola Norm.
Sup. Pisa} {\bf IV - 26} (1998), 133-160.
\v
\i{[BLF]} A.~Bressan and P.~LeFloch,  
Structural stability and regularity of entropy solutions to
hyperbolic systems of conservation laws, {\it Indiana Univ. Math. J.},
{\bf 48} (1999), 43-84.
\v
\i{[BG]} A.~Bressan and P.~Goatin,
Oleinik type estimates and uniqueness for $n\times n$
conservation laws, {\it J. Differential Equations} {\bf 156} (1999), 26-49.
\v
\i{[BLY]} A.~Bressan, T.~P.~Liu and T.~Yang,  $ L^1$ stability
estimates for $n\times n$ conservation laws, {\it Arch. Rational
Mech. Anal.} {\bf 149} (1999), 1-22. 
\v
\i{[BLY]} A.~Bressan, and T.~Yang,
On the convergence rate of vanishing viscosity approximations,
to appear.
\v
\i{[GL]} J.~Glimm and P.~Lax, Decay of solutions of systems of 
nonlinear hyperbolic 
conservation laws, {\it Amer. Math. Soc. Memoir} {\bf 101} (1970).
\v
\i{[L1]} T.~P.~Liu,
Decay to N-waves of solutions of general systems of nonlinear
hyperbolic conservation laws, {\it Comm. Pure Appl. Math.}
{\bf 30} (1977), 585-610.
\v
\i{[L2]} T.~P.~Liu, Admissible solutions of hyperbolic conservation
laws, {\it Amer. Math. Soc. Memoir} {\bf 240} (1981). 
\v
\i{[O]} O.~Oleinik, Discontinuous solutions of nonlinear
differential equations, {\it Amer. Math. Soc. Transl.} {\bf 26}
(1963), 95-172.
\bye